\documentclass[11pt,twoside,letter]{article}
\usepackage{amssymb}
\usepackage{hyperref}
\usepackage{color}
\usepackage{fancyhdr}

\newtheorem{thr}{Theorem}[section]
\newtheorem{lem}{\bf Lemma}[section]

\newtheorem{z}{Remark}[section]
\newcommand{\goth}{\Bbb}

\newcommand{\n}[1]{\refstepcounter{equation}\label{#1}
\eqno{(\arabic{section}.\arabic{equation})}}

 {}
\title{Controllability function as time of motion. I
\footnote{This is a translation of the paper: A.E. Choque Rivero,
V.I. Korobov, V.O. Skoryk: Controllability function as time of
motion I, (in Russian)  Mat. Fiz. Anal. Geom, 11(2), (2004),
208--225.}}
\author {\Large A.E. Choque-Rivero$^{\rm a}
\footnote{A.E. Choque-Rivero is currently professor and researcher
at the Instituto de F\'isica y Matem\'aticas, Universidad Michoacana
de San Nicol\'as de Hidalgo, Morelia, Mich., M\'exico,\newline
email: abdon@ifm.umich.mx},$ V.I. Korobov$^{\rm b,c}$,
 V.A. Skoryk$^{\rm b}$}

\begin{document}
 \maketitle
\center{\begin{tabular}{c}
 $^{\rm a}$Universidad Aut\'onoma del Carmen,\\
Calle 56 No. 4, C.P. 24180, Cd. del Carmen, Campeche, M\'exico\\
{\rm E-mail: achoque@pampano.unacar.mx, achoque@yahoo.com}\\
$^{\rm b}$Kharkov national university,\\
         Sq. Liberty, 4, 61077, Khakkov, Ukraine\\
       {\rm E-mail: $\{$vkorobov,skoryk$\}$@univer.kharkov.ua}\\
$^{\rm c}$Universitet Szczecinski, Institut Matematyki,\\
Wielkopolska str., 15, Institute of Mathematics, Szczecin, Poland\\
{\rm E-mail: korobow@sus.univ.szczecin.pl}
%%%%%%%%%%%%%%%%%%%%%%%%%%%%%%%%%%%
%%%%%%%%%%%%%%%%%%%%%%%%%%%%%%%%%%%-->   \'o ---- bez, Campeche
\end{tabular}}
%%%%%%%%%%%%%%%%%%%%%%%%%%%%%%%%%%%%%%%%%%%%%%%%%%%%%%%%%%%%%%%%%%%
%%%%%%%%%%%%%%%%%%%%%%%%%%%%%%%%%%%%%%%%%%%%%%%%%%%%%%%%%%%%%%%%%%
\center{{Mathematics Subject Classification 2000}: 93P50}\\
 \abstract{The
admissible positional control problem  for the canonical system with
geometrical restrictions on the control is considered.
 The investigation is performed with the help of the controllability
 function method. We obtain controllability functions which are
 the time of motion from an arbitrary initial condition to the
 origin.  We also reveal a set of controls which solves this problem.}
\thispagestyle{empty}
\newpage

\section{Introduction}
 The problem of synthesis with bounded controls (SBC) of the
  time optimal control \cite{pontryagin} is one of the
  well--known optimal control problems. The problem of SBC is stated
  as follows:

Given a system of differential equations,
$$
\dot{x} = f(x,u), \quad x \in {\goth R}^n,\quad u \in \Omega \subset
{\goth R}^r,\quad 0\in {\rm int}\ \Omega. \n{r1_l1}$$ It is required
to find a control $u$ of the form $ u=u(x),$ with values in the set
$ \Omega,$ such that the trajectory of the system
$$ \dot{x}=f(x,u(x)),\n{r1_l2}$$
starting at an arbitrary point $x_0,$ terminates at a given point
$x_1$ in minimal time.

\setcounter{page}{209}

In this case, the control $u(x)$ is time optimal and satisfies the
Bellman equation
 \cite{bellman}:
$$\min \limits_{u \in \ \Omega}\sum\limits_{i=1}^n\frac{\partial T(x)}{\partial x_i}f_i(x,u)=-1,\n{r1_l3}$$
where  $T(x)$ is the time of motion along the trajectory of  system
(\ref{r1_l2}),
 which corresponds to the control  $u(x).$  Let us denote by $\dot
T(x)_{|(\ref{r1_l1})},$ $\dot T(x)_{|(\ref{r1_l2})}$
 the time derivative along  system (\ref{r1_l1}), (\ref{r1_l2})
 respectively.  Hence equality  (\ref{r1_l3}) takes the form
$$\min \limits_{u \in \ \Omega}\dot T(x)_{|(\ref{r1_l1})}=\dot
T(x)_{|(\ref{r1_l2})}=-1,$$
 which means that the  derivative by virtue of  system
 (\ref{r1_l2}) of the time optimal $T(x)$ from an arbitrary point  $x$
 to a given point $x_1$ is equal to $-1.$
 Obviously this equality holds at points where the derivative
 exists.

 If we renounce the time optimality, we would consider the
 admissible positional synthesis problem, which consists of the
 construction of the positional control $u=u(x),$ which satisfies a
 given
  restriction, i.e.  $u \in \Omega$. As a result, %and such that
  the trajectory of
   system (\ref{r1_l2}), starting at an arbitrary point $x_0,$ terminates at
   a given point  $x_1$ at some {\it finite}  time  $T(x_0).$
 Furthermore, we assume that $x_1=0$ and $ f(0,0)=0.$

For solving this problem, the method of the controllability function
was introduced by V.I. Korobov \cite{korobov}. This method is based
in the construction of the controllability function
 $\;\Theta (x)$ $(\Theta (x)>0$ for $x\not=0$ and $\Theta
(0)=0)$, as well as, positional control $u(x)=\widetilde u
(x,\Theta(x)),$ such that  inequality
$$
\sum\limits_{i=1}^{n}\frac{\partial\Theta (x)}{\partial x_i}
f_i(x,u(x))\le -\beta\Theta^{1-1/\alpha}(x) \n{l3}$$ is satisfied
for some $\beta>0,$ $ \;\alpha>0.$
% If this inequality is fullfil
  The realization of this condition guarantees that the trajectory
  arrives to the origin at finite time.
 We assume that the controllability function
 $\Theta(x)$ is continuously differentiable at $x\ne 0.$
 In case that inequality (\ref{l3}) holds for
 $\alpha=\infty$, the function $\Theta(x)$ is a Lyapunov function.

In \cite{korobov_sklyar}
 a general method of constructing a controllability function and a
 control, which is solution of the admissible positional synthesis,
  was introduced for the linear control system
$$ \dot x=Ax+Bu,\quad x\in{\goth R}^n,\;\; u\in{\goth
R}^r,\n{l6}$$ where $A$ and $B$
 are $(n{\times}n)$ and $(n{\times}r)$ constant matrices,
 respectively.  In particular,
the function $\varphi (s)= 1{-}s$ for $ s\in[0,1],$ $\varphi (s)=0$
for $ s>1$ generates a controllability function
 $\Theta_\varphi (x),$ which is the time of  motion
 (\ref{l6}) with the control $u_\varphi(x)$ from an arbitrary point $x$
  to the origin according to the system (\ref{l6}).

We are interested in finding a wider set of pairs of functions: the
function of controllability $\Theta(x),$ which is the time of motion
and the control which solves the synthesis problem.

 In the first part of the work, we consider in detail this problem
 for the canonical system.
The second part is devoted to the solution of the positional
synthesis control problem for the completely controllable system
(\ref{l6}) with a restriction of the form $\;u\in \Omega=\{u:
\|u\|\le d\},$ $\;d>0.$
 The controllability functions $\Theta(x),$ which are the time
 of motion from the initial point $x$ to the origin,
 satisfy the equality $T(x)=\Theta(x),$
%then in this case
 and
the derivative of the controllability function by virtue of the
system $$\dot x=Ax+Bu(x)\n{l8}$$ is equal to $-1$:\,
 $\dot \Theta(x)_{|(\ref{l8})}=-1.$

 In the present work, we use the controllability function method
 which was first introduced in \cite{korobov}.
\section{Preliminary results}
\setcounter{equation}{0}
 We consider the canonical control system
$$ \dot{x}=A_0x+b_0u, \n{l10}$$ where $$ A_0= \left(
\begin{array}{cccccc}
 0    & 0    & \ldots  & 0 & 0 &0    \\
 1    & 0    & \ldots  & 0 & 0 & 0    \\
 \ldots & \ldots    & \ldots  &  \ldots&\ldots &\ldots  \\
 0   & 0 & \ldots     &0     & 1     &  0
\end{array}\right), \quad  b_0= \left(
\begin{array}{c}
 1     \\
 0   \\
 \ldots  \\
 0
\end{array}
\right),$$ with restrictions on    control $|u|\le d.$
   We will choose a control as a function depending on states coordinates
  $u=u(x)$ such that the trivial solution of  system (\ref{l10})
 for  $u=u(x)$ would be asymptotically stable. %This control
 We will refer to this control as the auxiliary
 control. Real numbers
 $a_1,a_2,\ldots,a_n$
 will be chosen in such a way that
${\lambda}^n-a_n{\lambda}^{n-1}-...-a_1=0$
 would have roots with the negative real part. In addition,
  set $u(x)=\sum\limits_{i=1}^n{a_i}x_i=(a,x),$ where
$a=(a_1,\ldots,a_n)^*.$ In this case,  system  (\ref{l10})
 has the form
$$\dot x =A_1x,\n{g1_2_f3}$$
where
  $$A_1=\left(\begin{array}{cccccc}
a_1 & a_2 &  \ldots & a_{n-2}& a_{n-1}& a_n\\
1 & 0 &  \ldots & 0& 0 & 0\\
\ldots &\ldots &\ldots &\ldots &\ldots&\ldots\\
0 & 0 &  \ldots &0 &1 &0
\end{array}\right);$$
its trivial solution is asymptotically stable.
 In this way, the auxiliary control solves the stabilizability
 problem for system (\ref{l10}) in all the phase space, but does not
 satisfy the given restriction.
  Since $A_1$ is stable, we can find a positive definite quadratic  form
$V(x)=(Fx,x)$ which is a Lyapunov function.
 The derivative of $V$ by virtue of  system
(\ref{g1_2_f3}) represents a negative definite quadratic form
$(-Wx,x),$ with positive definite matrices  $F$ and $W$.
  Since
$$\frac{d}{dt}(Fx,x)=(F\dot x,x)+(Fx,\dot x)=((FA_1+A_1^*F)x,x),$$
 this equation is reduced then to the Lyapunov matrix equation
$$FA_1+A_1^*F=-W,$$
 which for given matrices
 $A_1$ and $W$ has a  unique solution.

 Let  $F$ be a positive definite matrix, the form of which %form we
  will be specified later.
 Denote by $D(\Theta)$, $H$ diagonal matrices of the form
$$D(\Theta)={\rm diag}\left(\Theta^{-\frac{2i-1}{2}}\right)_{i=1}^n,\quad
H={\rm diag}\left(-\frac{2i-1}{2}\right)_{i=1}^n.$$

 Let $a_0$ be a positive number which will be determined later.
 For each $x\not=0$, define the controllability function $\Theta(x)$
 as a solution of the following equation:
$$2a_0\Theta=(D(\Theta)FD(\Theta)x,x).\n{g1_2_f8}$$
$$$$
 Set $\Theta(0)=0.$ This equality and equation
   (\ref{g1_2_f8}) determine the function $\Theta(x)$ which is
   continuous for
   all  $x$ and continuously differentiable for $x\ne0$.

We will look for control $u(x)$  in the form $$ u(x)=
\Theta^{-\frac{1}{2}}(x) a^* D(\Theta(x))x=\sum\limits_{k=1}^n\frac{
a_kx_k}{\Theta^k(x)},\n{r3_u}$$ where $F$ and the vector $a$ will be
 selected in a such a way that the controllability function
$\Theta(x)$ must  be the time of motion from the point
  $x$ to the origin. In other words, the time derivative $\dot \Theta(x)$
   by virtue of system (\ref{l10}) with the control
 $u(x)$  (\ref{r3_u}) satisfies equality
$$\dot\Theta(x)=-1.\n{dot_theta}$$
\begin{lem}   \label{le2.1}\hskip-2mm. %lemma2.1
 {\it Let the controllability function $\Theta(x)$
 satisfy equality (\ref{dot_theta}). Thus,
    $$a_1=-\frac{n(n+1)}{2} .\n{a_1}$$}
\end{lem}
 {\bf Proof}  Set $y(\Theta,x)=D(\Theta)x.$
 The controllability function  $\Theta(x)$ for
 $x\ne 0$ then satisfies equality
$$2a_0\Theta(x)=(Fy(\Theta(x),x),y(\Theta(x),x)),\n{ur_theta}$$
 and control
(\ref{r3_u}) has the form
$$u(x)=\Theta^{-\frac{1}{2}}(x)a^*y(\Theta(x),x).\n{r3_upr}$$
 In view of equalities $$ D(\Theta)A_0D^{-1}(\Theta)=\Theta^{-1}A_0,\quad
D(\Theta)b_0=\Theta^{-\frac{1}{2}}b_0$$
 the derivative of the function $y(\Theta(x),x)$  by virtue of the system
  (\ref{l10}) with this control has the form
$$\dot
y(\Theta(x),x)=\Theta^{-1}(x)\left(\dot\Theta(x)H+A_0+b_0a^*\right)y(\Theta(x),x).$$
 Therefore, from  equality  (\ref{ur_theta})  we find that
 the derivative of the controllability function by virtue of
 system (\ref{l10}) with control $u(x)$  (\ref{r3_upr})
 is given by equality
$$\dot \Theta(x)=\frac{\left((F(A_0+b_0a^*)+(A_0+b_0a^*)^*F)y(\Theta(x),x),y(\Theta(x),x)\right)}
{((F-HF-FH)y(\Theta(x),x),y(\Theta(x),x))}.$$

\newpage
  Hence, by the conditions of the lemma we have
   $$ F\left(A_0+b_0 a^*+\frac{1}{2}I-H\right)+
 \left(A_0+b_0a^*+\frac{1}{2}I-H\right)^*F=0,\n{l_1.10}$$
where $I$ is the identity matrix.  From this equality, we conclude
 that the matrix $F^\frac{1}{2}(A_0+b_0a^*+\frac{1}{2}I-H)F^{-\frac{1}{2}}$
 is a skew-symmetric matrix, %and
  therefore its characteristic
 equation has the form
$$ \det\left(F^\frac{1}{2}\left(A_0+b_0
a^*+\frac{1}{2}I-H\right)F^{-\frac{1}{2}}-\lambda I\right)= $$
$$=(\det F^{-\frac{1}{2}})^2\det\left(F\left(A_0+b_0
a^*+\frac{1}{2}I-H\right)-\lambda F\right)=0,
$$
 which is equivalent to
$$\det\left(\left(A_0+b_0 a^*+\frac{1}{2}I-H\right)-\lambda
I\right)=0. \n{l_1.80}
$$
In addition, it has eigenvalues with zero real part.  Let us write
equality (\ref{l_1.80}) in the following form:
$$ 0=\det \left(
\begin{array}{ccccc}
\lambda-1-a_1 & -a_2    & \ldots    & -a_{n-1} & -a_{n}    \\
- 1   &  \lambda-2   & \ldots  & 0 & 0     \\
\ldots & \ldots & \ldots    &\ldots   & \ldots \\
0 & 0  &\ldots & -1  & \lambda-n
\end{array}
\right)= $$
$$=\prod\limits_{j=1}^{n}(\lambda-j)-\sum\limits_{j=1}^{n-1}  a_j
\prod\limits_{i=j+1}^{n}(\lambda-i)-a_n=\lambda^n-\left(a_1+n(n+1)/2\right)\lambda^{n-1}-
\ldots.
$$
 From the fact that the roots of the last equality have zero real
 part and are complex conjugate, we have that
$a_1=-n(n+1)/2.$

%lem1.2
\begin{lem}\label{le2.2}\hskip-2mm. %lemma2.2
{\it Let  $F=\left(f_{ij}\right)_{i,j=1}^n$ be a positive definite
matrix.

Thus, the matrix
$$G(\xi)=\left (\begin{array}{cccc}
\xi & f_{12}&\ldots & f_{1n}\\
f_{21}& f_{22} & \ldots & f_{2n}\\
\ldots &\ldots & \ldots & \ldots\\
f_{n1} & f_{n2}&\ldots & f_{nn}
\end{array}\right)%\n{m_g}
$$
%with $\xi\ge f_{11}$
 is a positive definite matrix and $\det G(\xi_0)=0$ for $\xi>\xi_0$.}
\end{lem}
\newpage
 {\bf Proof.} The determinant $\det\, G(\xi)$ is a linear function on
 $\xi$, i.e., $G(\xi)=a\xi+b$, hence by Sylvester criteria
 $a>0$. Therefore, if $\det\, G(f_{11})=a\,f_{11}+b>0$, then for
  any $\xi\geq f_{11}$ we obtain that $\det\, G(\xi)$. Moreover,
  there exists a unique $\xi_0$ such that
  $\det\,G(\xi_0)=a\,\xi_0+b=0$,
  and obviously, for any $\xi>\xi_0$ we have that $\det\,G(\xi)>0$.

 \begin{lem}\label{le2.3}\hskip-2mm. %lemma 2.3
 {\it
 \begin{itemize}
\item[i)] The $(s{\times}s)$ matrices
$$ P_{s,k}=\left(\frac{1}{i+j+k-2}\right)_{i,j=1}^s= \left(
\begin{array}{ccccc}
\frac{1}{k} & \frac{1}{k+1} & \ldots & \frac{1}{s+k-1}\\
\frac{1}{k+1}& \frac{1}{k+2} &\ldots & \frac{1}{s+k}\\
\ldots& \ldots           & \ldots & \ldots\\
\frac{1}{s+k-1}     &\frac{1}{s+k} &\ldots & \frac{1}{2s+k-2}
\end{array}
\right), \quad k\in{\goth N}, $$ are positive definite ones.
\item[ii)] For every $k\in {\goth N}$, we have $$\Delta_{s,k}={\rm det}
\left(\frac{1}{i+j+k-2}\cdot \frac{1}{i+j+k-1}\right)_{i,j=1}^s>0.$$
\item[iii)] The determinant of the  $(n{\times}n)$ matrix $P$ of the form
$$ P= \left(
\begin{array}{ccccc}
1-\frac{n(n+1)}{2} &\frac{1}{2\cdot 3}& \ldots & \frac{1}{n(n+1)}\\
-1    &\frac{1}{3\cdot 4}& \ldots & \frac{1}{(n+1) (n+2)}\\
0     &\frac{1}{4\cdot 5}& \ldots & \frac{1}{(n+2)(n+3)}\\
\ldots& \ldots           & \ldots & \ldots\\
0     &\frac{1}{(n+1) (n+2)} &\ldots & \frac{1}{(2n-1) 2n}
\end{array}
\right). \n{l1_13}$$ is equal to zero, and its rank is equal to
 $(n{-}1).$
\end{itemize}}
\end{lem}
{\bf Proof.} i) The positive definiteness of the matrices
 $P_{s,k}$ follows from the positive definiteness of the
  Hilbert matrix
$\left(\frac{1}{i+j-1}\right)_{i,j=1}^s.$

 The proof of this fact readily follows from the integral
 representation
$$
P_{s,k}=\int\limits_0^1 \left(
\begin{array}{c}
1\\t\\ \ldots\\ t^{s-1}
\end{array}
\right) (1, t, \ldots, t^{s-1})d\frac{t^k}{k}.
$$
The matrix $P_{s,k}$ is then positive definite one
 \cite{krein}.

ii) Let us represent the determinant $\Delta_{s,k} $ in the form
$$\Delta_{s,k}={\rm det}
\left(\frac{1}{i+j+k-2}- \frac{1}{i+j+k-1}\right)_{i,j=1}^s.$$
\newpage
 Add to the  $i$-th row
($i\ge 1$) all the next rows, that is, the rows with numbers
 $i{+}1,$ $\ldots,$ $s,$ and extract
   from every
$i$-nth row  ($i=1,\ldots,s$) of the obtained determinant
 the value
  $(s{-}i{+}1)$ and from every $j$-th column the
  value $1/(s{+}j{+}k{-}1)$ $j=1,\ldots,s.$
 We obtain
$$\Delta_{s,k}= \frac{s!}{(s+k)\ldots(2s+k-1)} \det P_{s,k}.$$
From there and in view of i), it follows that $\Delta_{s,k}>0,$
$k\in{\goth N}.$

 iii)  Let us write the determinant of the matrix $P$ in the form $$
\left|
\begin{array}{ccccc}
1-\frac{n(n+1)}{2} &\frac{1}{2}-\frac{1}{3}& \ldots &
\frac{1}{n}-\frac{1}{n+1}\\
-1    &\frac{1}{3}-\frac{1}{4}& \ldots & \frac{1}{n+1}-\frac{1}{n+2}\\
\ldots& \ldots           & \ldots & \ldots\\
0     &\frac{1}{(n+1)}-\frac{1}{n+2} &\ldots &
\frac{1}{2n-1}-\frac{1}{2n}
\end{array}
\right|. \n{l1_14}$$
 We transform this determinant in two different ways.
 In the beginning, we add to every $j$-th column ($j\ge 2$)
 of this determinant all the next ones;
 the columns with numbers $j{+}1,$
$\ldots,$ $n.$ We have $$ {\det} P= \left|
\begin{array}{ccccc}
1-\frac{n(n+1)}{2}&\frac{n-1}{2(n+1)}&\ldots &\frac{2}{(n-1)(n+1)}&
        \frac{1}{n(n+1)}\\
-1    &\frac{n-1}{3(n+2)}&\ldots &\frac{2}{n(n+2)}&\frac{1}{(n+1)(n+2)}\\
\ldots     &\ldots              & \ldots  & \ldots & \ldots\\
0     &\frac{n-1}{(n+1)2n} &\ldots &\frac{2}{(2n-2)2n}&
\frac{1}{(2n-1)2n}
\end{array}\right|=$$
$$=\frac{(n-1)!}{(n+1)\ldots 2n} \left|
\begin{array}{ccccc}
(n+1) \left(1-\frac{n(n+1)}{2}\right)
  &\frac{1}{2}& \ldots &\frac{1}{n-1}& \frac{1}{n}\\
-(n+2)    &\frac{1}{3}&\ldots &\frac{1}{n}&\frac{1}{n+1}\\
\ldots &\ldots& \ldots& \ldots            & \ldots\\
0     &\frac{1}{n+1} &\ldots &\frac{1}{2n-2}& \frac{1}{2n-1}
\end{array}\right|. \n{l1_15}$$
 On the other side, by adding the  rows with numbers
 from $i+1$ till  $n$
 %next numbers
 to the $i$-th row
  ($i\ge 1$) of determinant (\ref{l1_14}), we obtain
 % hhhhhhhhhhhhhhhhhhhhhhhhhhhhhhhhhhhhhh
 % hhhhhhhhhhhhhhhhhhhhhhhhhhhhhhhh
 %of  determinant (\ref{l1_14}), the next ones,
 %that is, the rows with numbers form $i+1$ till  $n$, we obtain
$${\det P= \left|
\begin{array}{cccccc}
-\frac{n(n+1)}{2}   &\frac{n}{2(n+2)}&
\cdots &\frac{n}{(n-1)2n}& \frac{n}{n(2n)}\\
-1    &\frac{n-1}{3(n+2)}&\ldots
 &\frac{n-1}{n(2n)}&\frac{n-1}{(n+1)2n}\\
\ldots &\ldots& \ldots& \ldots & \ldots\\
0     &\frac{1}{(n+1)(n+2)}  &\ldots &\frac{1}{(2n-2)(2n-1)}&
\frac{1}{(2n-1)2n}
\end{array}\right|}=$$
$$=\frac{n!}{(n+2)\ldots2n} \left|
\begin{array}{cccccc}
-\frac{n+1}{2}
  &\frac{1}{2}& \ldots &\frac{1}{n-1}& \frac{1}{n}\\
-\frac{1}{n-1}    &\frac{1}{3}&\ldots &\frac{1}{n}&\frac{1}{n+1}\\
0 &\ldots& \ldots& \ldots            & \ldots\\
0     &\frac{1}{n+1} &\ldots &\frac{1}{2n-2}& \frac{1}{2n-1}
\end{array}\right|. \n{l1_16}$$
\newpage
 Denote by $B_1$ and $B_2$
 the complementary minors with numbers (1,1) and (1,2),
 respectively, which appear in the RHS of
(\ref{l1_15}) and (\ref{l1_16}). Thus
$$\det P=\left(\left(1-\frac{n(n+1)}{2}\right)\  B_1+\frac{n+2}{n+1}\ B_2\right)\frac{(n-1)!}{(n+1)\cdots 2n},$$
$$\det P=\left(-\frac{n(n+1)}{2}\ B_1+\frac{n}{n-1}\ B_2\right)\frac{(n-1)!}{(n+1)\cdots 2n},$$
and
$$
 B_1=\frac{2}{n^2-1}\ B_2,\quad \det P=0.$$
 Since
   $ B_1\ne 0,$  then  ${\rm rank} (P)=n-1.$
 The Lemma is proved.

 We determine  matrix $F,$ which appears in equation (\ref{g1_2_f8}),
 and the vector $a,$ which appears in control of  form
  (\ref{r3_u}). From equality  (\ref{l_1.10}), we have
$$
\left(A_0+b_0 a^*+\frac{1}{2}I-H\right)F^{-1} +F^{-1}\left(A_0+b_0
a^*+\frac{1}{2}I-H\right)^*=0. \n{l1_18}$$

 Denote
 $D_n={\rm diag} \left((-1)^{i-1}/{(i-1)!}\right)_{i=1}^{n}.$
  We must find
  the matrix  $F^{-1}$ with the form  $F^{-1}=D_nCD_n,$
 where the matrix $C$ is a Hankel matrix
  $C=\left(c_{i+j}\right)^{n-1}_{i,j=0}.$
 We find the matrix  $C.$ To this end
 we write equality (\ref{l1_18}) in the form
$$\left(\frac{1}{2}I-H{+}D^{-1}_nA_0D_n+D_n^{-1}b_0a^*D_n\right)C+$$
$$+C\left(\frac{1}{2}I-H{+}D^{-1}_nAD_n+D_n^{-1}b_0 a^*D_n \right)^*=0.  \n{l1_19}
$$ Denote by $Q=\left(q_{ij}\right)_{i,j=1}^n$ the matrix in the LHS of
 equality (\ref{l1_19}) and write it in the following form
$$
Q= \left(\begin{array}{cc} Q_{11} & Q_{12}\\Q_{21} & Q_{22}
\end{array}\right),$$
 where
$Q_{11}=q_{11},$ $ Q_{12}=(q_{12},\ldots, q_{1n}),$ $
Q_{21}=Q_{12}^*,$ $Q_{22}=\left(q_{ij}\right)_{i,j=2}^n.$
 Hence the elements of the matrix $Q$ are given by
the following equality:
$$\left\{\begin{array}{l}
\displaystyle q_{11} = 2(1+\widetilde a_1)c_0+2\sum\limits_{j=1}^{n-1}\widetilde a_{j+1}c_j, \\[10pt]
\displaystyle q_{1i} =-(i-1)c_{i-2}+(1+\widetilde a_1+i)c_{i-1}+
\sum\limits_{j=1}^{n-1}\widetilde a_{j+1}c_{i+j-1},\;
 i=2,\ldots,n,\\[10pt]
\displaystyle q_{ij} = -(i+j-2)c_{i+j-3}+(i+j)c_{i+j-2},\;
i=2,\ldots,n,\; j=2,\ldots,n,
\end{array}\right. \n{x2}$$
\newpage
 where  $\widetilde a_1,$ $\ldots,$ $\widetilde a_n$
 are the components of the vector
  $\widetilde a=D_na.$
 Since $Q$ is a zero matrix, then from the equalities $Q_{22}=0$
 we obtain
$$c_j=\frac{(2n-1)2n}{(j+1)(j+2)}c_{2n-2},\quad j=1,\ldots,2n-3,$$
 and equality (\ref{x2}) has de form
$$\begin{array}{l}
\displaystyle \frac{(1+\widetilde
a_1)c_0}{2n(2n-1)c_{2n-2}}+\sum\limits_{j=1}^{n-1}
     \frac{\widetilde a_{j+1}}{(j+1)(j+2)}=0,  \\[10pt]
\displaystyle -\frac{c_0}{2n(2n-1)c_{2n-2}}+\frac{3+\widetilde
a_1}{2\cdot 3}
+\sum\limits_{j=1}^{n-1} \frac{\widetilde a_{j+1}}{(j+2)(j+3)}=0,\\[10pt]
\displaystyle -\frac{1}{3}+\frac{4+a_1}{3\cdot
4}+\sum\limits_{j=1}^{n-1}
\frac{\widetilde a_{j+1}}{(j+3)(j+4)}=0,\\[10pt]
.\qquad.\qquad.\qquad.\qquad.\qquad.\qquad.\qquad.\qquad.\qquad.\qquad.\\
 \displaystyle
-\frac{1}{n}+\frac{n+1+\widetilde
a_1}{n(n+1)}+\sum\limits_{j=1}^{n-1} \frac{\widetilde
a_{j+1}}{(j+n)(j+n+1)}= 0, \end{array} \n {system}$$
 where $\widetilde a_1=a_1$ is defined by equality (\ref{a_1}).
 Consider equality  (\ref{system}), as a linear
 system of equation with the respect to the unknowns
 % как линейную систему уравнений
%относительно неизвестных
$$\frac{c_0}{(2n-1)2n c_{2n-2}},\quad \widetilde a_2,\quad \ldots,\quad
 \widetilde a_n,$$
 which in vector form can be written as
%которая в векторной форме записи имеет  вид
$$Py=y_0,\n{slau}$$ where  $P$  is a $(n{\times}n)$ matrix defined
by equality (\ref{l1_13}). In addition,   $y,$ $y_0$ are vectors of
the form
$$y{=}\left(\frac{c_0}{(2n-1)2n c_{2n-2}}, \widetilde a_2,\ldots, \widetilde a_n
\right)^*,\; y_0{=}\left(0,-\frac{3+a_1}{2\cdot 3},
-\frac{a_1}{3\cdot 4}\ldots,
     -\frac{a_1}{n(n+1)}\right)^*.$$
     This system represents a subsystem, since the matrix
 $C=\left(
\frac{1}{(i+j-1)(i+j)}\right)_{i,j=1}^{n}$
 and the vector
$a^*=-\frac{1}{2}b_0^*D_n^{-1} C^{-1} D_n^{-1}$
 satisfy the equality
  (\ref{l1_19}). Actually, the equality  (\ref{l1_19})
 for such a
%%%%%%%%%%%%%%%%%%%%%%%%%%%%%%%%%%%%%%%%
  selection of the vector $a$  has the form
$$\left(\frac{1}{2}I-H+D^{-1}_nA_0D_n\right)C+
C\left(\frac{1}{2}I-H+D^{-1}_nA_0D_n \right)^*=b_0b_0^*.$$
 This equality holds if one inserts the matrix $C$
  in the last equality.

 Therefore,
 $ {\rm rank}\,  (P)= {\rm rank} (P, y_0)$;
  by virtue of  iii) and Lemma  \ref{le2.3} we have that
 ${\rm rank}\, (P)=n-1.$
\newpage
 Let us find all the solutions of system
  (\ref{slau}). Consider matrix  $((n{-}1){\times}(n{-}1))$,
  which has the form $\widetilde P$
$$\widetilde P=
\left(
\begin{array}{cccc}
-1    &\frac{1}{3 \cdot 4}& \ldots & \frac{1}{n (n+1)}\\
0     &\frac{1}{4 \cdot 5}& \ldots & \frac{1}{(n+1) (n+2)}\\
\ldots& \ldots           & \ldots & \ldots \\
0     &\frac{1}{(n+1)(n+2)} &\ldots & \frac{1}{(2n-2)(2n-1)}
\end{array}\right).$$
 This matrix, by because of ii) of Lemma \ref{le2.3}
 for  $k=4$ and $s=n{-}2,$ is a nondegenerate matrix.
 Denote by $\widetilde y_0= a_n d'+d'',$ where $$
 d'= \left(
-\frac{1}{(n+1)(n+2)},\; -\frac{1}{(n+2)(n+3)},\; \ldots,\;
     -\frac{1}{(2n-1)2n}
\right)^*, $$ $$ d''= \left( -\frac{3+a_1}{2\cdot
3},\;-\frac{a_1}{3\cdot 4},\; \ldots,
   \;  -\frac{a_1}{n(n+1)}\right)^*.
$$
Furthermore,
 by  $\widetilde y$ we denote the  $(n{-}1)$-dimensional vector
 $$ \widetilde y=
\left( \frac{c_0}{(2n-1)2n c_{2n-2}},\;\widetilde a_2,\; \ldots,\;
 \widetilde a_{n-1}
\right)^*$$
 and consider a system of equations $\widetilde P
\widetilde y=\widetilde y_0 $ with respect to
  $\widetilde y$.
 This system has a unique solution
   $\widetilde
y=\widetilde P^{-1}\widetilde y_0$:
$$\begin{array}{l}
\displaystyle  \frac{c_0}{(2n-1)2n\, c_{2n-2}}
=\frac{1}{\Delta}\left(\Delta'_1 \widetilde a_n + \Delta''_1
\right)=\frac{1}{\Delta}\left(\Delta'_1
\frac{(-1)^{n-1}}{(n-1)!} a_n + \Delta''_1 \right),\\[15pt]
\displaystyle \widetilde  a_j= \frac{1}{\Delta} \left(\Delta'_j
\widetilde a_n + \Delta''_j \right)=\frac{1}{\Delta} \left(\Delta'_j
\frac{(-1)^{n-1}}{(n-1)!}a_n + \Delta''_j \right),\quad
 i=2,\ldots, n{-}1, \end{array}            \n{l1_25}$$
where  $\Delta=\det\, \widetilde P,$ and $\Delta'_j,$ $\Delta''_j$
$(j=1,\ldots,n{-}1)$ are the determinants of the matrix
 $\widetilde P,$
 in which instead of its the $j$-th the columns, vectors
 $d',$ $d''$ are inserted. By virtue of condition
  ${\rm rank }\, (P)={\rm rank }\, (P,\;y_0)=n{-}1$
  and relations
 (\ref{l1_25}) describe all the solutions of the system
(\ref{slau}).

 Therefore, the next lemma follows.
\begin{lem}  \label{le2.4}\hskip-2mm. {\it  %lemma2.4
 Matrix $$ C{=}\left(
\begin{array}{cccc}
\frac{1}{\Delta}\left(\Delta'_1 \frac{({-}1)^{n{-}1}}{(n{-}1)!}a_n
{+}\Delta''_1\right) &\frac{1}{2\cdot 3}&\ldots&
   \frac{1}{n(n{+}1)}\\
\frac{1}{2{\cdot} 3} & \frac{1}{3{\cdot} 4} & \cdots&\frac{1}{(n{+}1)(n{+}2)}\\
\ldots & \ldots & \ldots& \ldots\\
\frac{1}{n(n{+}1)} &
\frac{1}{(n{+}1)(n{+}2)}&\ldots&\frac{1}{(2n{-}1)2n}
\end{array}
\right)(2n{-}1)2n\, c_{2n{-}2} \n{l1_26} $$
\newpage
 and the vector $$  a=\left(
-\frac{n(n+1)}{2},\;-\frac{1}{\Delta}\left(\Delta'_2
\frac{(-1)^{n-1}}{(n-1)!}a_n+\Delta''_2\right),\; \ldots,\;\right.$$
$$\left. \ldots,(-1)^{n-2}(n-2)! \frac{1}{\Delta}\left(\Delta'_{n-1}
\frac{(-1)^{n-1}}{(n-1)!}a_n+\Delta''_{n-1}\right),\;a_n\right)^*
\n{l1_27}$$
 give all the solutions of equation
  (\ref{l1_19}); moreover for
$$c_{2n-2}>0,\quad \frac{1}{\Delta}\left(\Delta'_1 \frac{(-1)^{n-1}}{(n-1)!}a_n +\Delta''_1\right)>
\xi_0,\n{usl_c_an}$$
 where
 $\xi_0$ is a root of the equation
$$\left|
\begin{array}{cccc}
\xi_0 &\frac{1}{2\cdot 3}&\cdots&
   \frac{1}{n(n+1)}\\
\frac{1}{2\cdot 3} & \frac{1}{3\cdot 4} & \cdots&\frac{1}{(n+1)(n+2)}\\
\cdots & \cdots & \cdots& \cdots\\
\frac{1}{n(n+1)} & \frac{1}{(n+1)(n+2)}&\cdots&\frac{1}{(2n-1)2n}
\end{array}\right|=0,\n{u_xi_0}$$
 matrix $C$ is a positive definite one.
 }
\end{lem}
{\bf Proof.} Let  $c_{2n-2}>0$, and consider  matrix
$$ \widetilde C(z)=\left(
\begin{array}{cccc}
z &\frac{1}{2\cdot 3}&\cdots&
   \frac{1}{n(n+1)}\\
\frac{1}{2\cdot 3} & \frac{1}{3\cdot 4} & \cdots&\frac{1}{(n+1)(n+2)}\\
\cdots & \cdots & \cdots& \cdots\\
\frac{1}{n(n+1)} & \frac{1}{(n+1)(n+2)}&\cdots&\frac{1}{(2n-1)2n}
\end{array}
\right)(2n-1)2n\, c_{2n-2}. $$
 This matrix is positive definite for
$z=1/2,$ by virtue of  ii) of the Lemma \ref{le2.3}, its principal
 minors  $\Delta_{s1},$ $s=1,\ldots,n,$ are positive:
the Sylvester criterion holds. Because of Lemma \ref{le2.2}, matrix
$\widetilde C(z)$  is positive definite for $z> \xi_0,$ where
   $\xi_0$ is a root of the equation
$\det\ \widetilde C(\xi_0)=0.$ From where we obtain that for the
parameters $c_{2n-2}$ и $a_n,$ which satisfy condition
 (\ref{usl_c_an}), matrix  $C$ is positive definite.

\begin{z}\label{r1_z1}\hskip-2mm. Observe that
equality
 (\ref{u_xi_0}) is reduced to the form
$$\left|
\begin{array}{cccc}
\left(1+\frac{1}{n}\right)\xi_0+\frac{1}{2}-\frac{1}{2n}
&\frac{1}{2}&\ldots&
   \frac{1}{n}\\
\frac{1}{2} & \frac{1}{3} & \cdots&\frac{1}{n+1}\\
\ldots & \ldots & \ldots& \ldots\\
\frac{1}{n} & \frac{1}{n+1}&\cdots&\frac{1}{2n-1}
\end{array}\right|=0.\n{u_xi_0_1}$$
 For establishing this fact, to each $i$-th row
 we add all the remaining rows;  from the
 obtained  rows and columns of the new determinant we then extract
   common factors. For the cases
$n=2,3,4,5,6$ and $7$ the roots of $\xi_0$ is equal to $1/3,$
$5/12,$ $9/20,$ $7/15,$ $10/21,$  and $27/56$ respectively.
\end{z}
\newpage

\begin{lem}  \label{le2.5}\hskip-2mm. {\it Matrix $C^1=C-CH-HC$
for the form
$$ C^1=\left(
\begin{array}{cccc}
\frac{2}{\Delta}\left(\Delta'_1 \frac{(-1)^{n-1}}{(n-1)!}a_n
+\Delta''_1\right) &\frac{1}{2}&\ldots&
   \frac{1}{n}\\
\frac{1}{2} & \frac{1}{3} & \ldots&\frac{1}{n+1}\\
\ldots & \ldots & \ldots& \ldots\\
\frac{1}{n} & \frac{1}{n+1}&\ldots&\frac{1}{2n-1}
\end{array}
\right)(2n-1)2n\, c_{2n-2} \n{lem_f0}$$
 is positive definite for
$$c_{2n-2}>0,\quad \frac{1}{\Delta}\left(\Delta'_1 \frac{(-1)^{n-1}}{(n-1)!}a_n
+\Delta''_1\right)>\left(\frac{1}{2}+\frac{1}{2n}\right)\xi_0
+\frac{1}{4}-\frac{1}{4n},$$
where $\xi_0$ is a root of equation (\ref{u_xi_0_1}).}
\end{lem}
{\bf Proof.} By  i) of Lemma \ref{le2.3} for  $s=n$ and $k=1$,
  matrix $C^1$ is positive definite if the parameter   $a_n$
   satisfies the condition
$$\frac{2}{\Delta}\left(\Delta'_1 \frac{(-1)^{n-1}}{(n-1)!}a_n +\Delta''_1\right)=1.$$
 From Lemma  \ref{le2.2} for $c_{2n-2}>0$, matrix  $$\widetilde
C^1(z)=\left(\begin{array}{cccc} z &\frac{1}{2}&\ldots&
   \frac{1}{n}\\
\frac{1}{2} & \frac{1}{3} & \ldots&\frac{1}{n+1}\\
\ldots & \ldots & \ldots& \ldots\\
\frac{1}{n} & \frac{1}{n+1}&\ldots&\frac{1}{2n-1}
\end{array}
\right)(2n-1)2n\, c_{2n-2} \n{lem_f2} $$
  is positive definite for all
$$z>\xi,\n{lem_f3}$$ where  $\xi$
is a root of equality $\det \ \widetilde C^1(\xi)=0.$
 By equality (\ref{u_xi_0_1}), we have
$$0=\det \ \widetilde C^1(\xi)= \det \ \widetilde
C^1\left(\left(1+\frac{1}{n}\right)\xi_0+\frac{1}{2}-\frac{1}{2n}\right),$$
  and consequently
$$\xi=\left(1+\frac{1}{n}\right)\xi_0+\frac{1}{2}-\frac{1}{2n}.\n{lem_f4}$$
 Thus, the proof readily follows by comparing (\ref{lem_f0})
  and (\ref{lem_f2}), as well as (\ref{lem_f4})
  and the condition (\ref{lem_f3}).
\newpage
\section{Synthesis of the bounded controls for the canonical system}
\setcounter{equation}{0}
 The solution of the synthesis problem for the canonical system
 (\ref{l10})  when the controllability function
 is  time of motion, gives the following solution

\begin{thr}  \label{thr1}\hskip-2mm.
 {\it
 Let  $$c_{2n-2}>0,\;
 \frac{1}{\Delta}\left(\Delta'_1
 \frac{(-1)^{n{-}1}}{(n{-}1)!}a_n {+}\Delta''_1\right){>}\max\left\{\xi_0,\;
 \left(\frac{1}{2}{+}\frac{1}{2n}\right)\xi_0{+}
 \frac{1}{4}{-}\frac{1}{4n}\right\},\n{r3_usl_par}$$
 where $\xi_0$  is a root of   equation
 (\ref{u_xi_0_1}).
 Let the matrix $C$ and the vector $a$ be as in  (\ref{l1_26})
  and (\ref{l1_27}), respectively.
  The number  $a_0$ satisfies the condition
 $$0<a_0\le \frac{d^2}{2(F^{-1}a,a)}; \quad F^{-1}=D_nCD_n,\n{a0}$$
 and the controllability function
  $\Theta(x)$ for $x\ne 0$ is defined by equality
    (\ref{g1_2_f8}),
 where $F=D_n^{-1}C^{-1}D_n^{-1},$ and $\Theta(0)=0.$

 Thus, control $u(x)$  (\ref{r3_u})
 transfers the arbitrary initial point
   $x\in
{\goth R}^n$ to the origin along the trajectory
 $\dot x= A_0x+b_0u(x)$ in time  $T(x)=\Theta(x)$
  and satisfies the restriction $|u(x)|\le d.$}
\end{thr}
{\bf Proof.}  By Lemma \ref{le2.4} and  condition
(\ref{r3_usl_par}), matrix $F^{-1}=D_nCD_n,$ and consequently matrix
$F$
 are positive definite.

  If matrix  $(F{-}FH{-}HF) $ is positive definite,
   equation (\ref{g1_2_f8})  for  $x\ne 0$ has a unique continuously
differentiable solution $\Theta(x)$; see   \cite{korobov}.

 We establish for which parameters $a_n$ and $c_{2n-2}$ is matrix
$(F{-}FH{-}HF)$ positive definite.
 The positive definiteness of this matrix
  will follow from the positive definiteness
   of matrix $(F^{-1}{-}HF^{-1}{-}F^{-1}H).$
 Since matrix
 $F^{-1}=D_nCD_n,$ and  matrices $D_n$ and $H$
  commute with each other, then the following equality is valid
$$F^{-1}{-}HF^{-1}{-}F^{-1}H=D_n(C-HC-CH)D_n.$$
By Lemma \ref{le2.5} and  condition (\ref{r3_usl_par}),  matrix
$(C{-}HC{-}CH)$ is positive definite;
 hence  matrix $(F{-}FH{-}HF)$ is also positive definite.

 Therefore, under conditions
 (\ref{r3_usl_par})  equation  (\ref{g1_2_f8}) for
$x\ne 0$ has a unique continuously differentiable solution
 $\Theta(x).$ Set
 $\Theta(0)=0,$ and $\Theta(x)$ becomes a continuous function for all  $x.$
\newpage
 Now let us establish that the control is bounded. We  estimate
  the expression
$a^*y(\Theta,x){\Theta}^{-\frac{1}{2}}.$
  To this end, we fix   $\Theta$
 and solve the problem of finding the extremal
 of the function
 $a^*y(\Theta,x){\Theta}^{-\frac{1}{2}}$ under the condition
$$(Fy(\Theta,x),y(\Theta,x))-2a_0\Theta=0.\n{g1_p2_ogr}$$
 This problem we solve with the help of the Lagrange multipliers.
 The Lagrange function has the form
$$a^*y(\Theta,x){\Theta}^{-\frac{1}{2}}-\lambda(Fy(\Theta,x),y(\Theta,x))+
2\lambda a_0\Theta.$$
 Let $y_0$ be the extremal point. The necessary condition
 gives
 $a \Theta^{-\frac{1}{2}}-2\lambda Fy_0=0,$ from where we have
  $\;\displaystyle
y_0=1/(2\lambda)\Theta^{-\frac{1}{2}}F^{-1}a.\;$
 Setting $y_0$ in restriction
(\ref{g1_p2_ogr}), we have $1/(2\lambda)=\pm \sqrt{2a_0/(F^{-1}a,
a)}\;\Theta.$ Consequently, we obtain
$(a,y_0){\Theta}^{-\frac{1}{2}}=\pm\sqrt{2a_0(F^{-1}a,a)}.$
 Hence, by using the control $u(x),$ we obtain
$$|u(x)|\le \sqrt{2a_0(F^{-1}a,a)}.\n{g1_2_g20_1}$$
 Selecting the number $a_0$ from condition (\ref{a0}),
 from  inequalities (\ref{g1_2_g20_1}) it follows that
 $u(x)$ satisfies the given restriction in all the phase space.

Hence from Theorem 1 \cite{korobov}, control $u(x)$
 solves the synthesis problem of bounded controls
 in all the phase space, and the time of motion from
  point $x$ to the origin is equal to the function of controllability
 at $x$: $T(x)=\Theta(x).$ \vskip5mm

 {\bf Example.} Let us consider the
  synthesis problem for the system
$$\dot x_1=u,\quad \dot x_2=x_1,\quad \dot x_3=x_2,\n{pr_l1}$$
 with a restriction on the control of the form
  $\;|u|\le 1.$

 Since
  $\Delta=-1/20,$ $\Delta_1'=1/3600,$ and
$\Delta_1''=-1/60, $ then vector  $a$  (\ref{l1_27}), matrix $C$
(\ref{l1_26}) and matrix  $D_3$ take the form
$$a=\left(\begin{array}{c}
-6 \\\frac{a_3}{3}-10\\ a_3\end{array}\right),\quad
C=\left(\begin{array}{ccc}
10-\frac{a_3}{12} & 5 &\frac{5}{2}\\
5 &\frac{5}{2} &\frac{3}{2}\\
\frac{5}{2} &\frac{3}{2}& 1
\end{array}\right)c_4,\quad D_3=\left(\begin{array}{ccc} 1 & 0 & 0\\ 0 & -1 & 0 \\ 0 &
0 & \frac{1}{2}
\end{array}\right).$$
Thus, the matrix $F^{-1}=D_3CD_3$ and its inverse matrix has the
form
$$F^{-1}{=}\left(\begin{array}{crr}
10-\frac{a_3}{12} & - 5 &\frac{5}{4}\\
-5 &\frac{5}{2} &-\frac{3}{4}\\
\frac{5}{4} &-\frac{3}{4}& \frac{1}{4}
\end{array}\right)c_4,\;
F{=}-\frac{c_4}{30{+}a_3}\left(\begin{array}{ccc}
12 & 60 &120\\
60 &180{-}4a_3 &240{-}12a_3\\
120 &240{-}12a_3& -40a_3
\end{array}\right).$$
\newpage
 From (\ref{u_xi_0}) and we obtain that $\xi_0=5/12.$
 Condition (\ref{r3_usl_par}) then has the form
$$c_4>0,\quad \frac{1}{3}-\frac{a_3}{360}>\frac{4}{9},$$
 according to it, we set
 $c_4=1,$ $a_3=-45.$ As a result,
$$a=\left(\begin{array}{c}
-6 \\-25\\ -45\end{array}\right),\quad
F^{-1}=\frac{1}{4}\left(\begin{array}{rrr}
55 & -20 &5\\
-20 & 10 &-3\\
5 &-3& 1
\end{array}\right),\quad
F=4\left(\begin{array}{rrr} 1/5 & 1 & 2\\ 1 & 6 & 13 \\ 2 & 13 & 30
\end{array}\right).$$
 We select the number  $a_0$ of  equality
   (\ref{g1_2_f8}) from conditions  (\ref{a0}), which in this
   case has the form $\; 0<a_0\le 2/205.\;$
 We determine the controllability function   $\Theta(x)$ for    $x\ne 0$
    as a positive solution of equation
 (\ref{g1_2_f8})  (which is unique).
 This equation in this case has the form
$$\Theta^6=41  \Theta^4x_1^2{+}410  \Theta^3x_1x_2{+}820 \Theta^2x_1x_3{+}
1230  \Theta^2 x_2^2{+} 5330  \Theta x_2x_3{+} 6150x_3^2.\n{l2_41}$$
 The control
 $u(x)$  (\ref{r3_u}), which solves the global synthesis problem
 for system
(\ref{pr_l1}) and satisfies the restriction
 $|u(x)|\le 1,$ is given by the formula
$$u(x)=-\frac{6}{\Theta(x)}x_1-\frac{25}{\Theta^2(x)}x_2-
\frac{45}{\Theta^3(x)}x_3. \n{l2_42} $$

Let us find the trajectory of the system (\ref{pr_l1}), which
corresponds to the control
 $u=u(x)$ of form  (\ref{l2_42}) and starts at the point
$x(0)=x_0\in{\goth R}^3$. This trajectory
 represents a solution of the system
$$ \dot x_1=-\frac{6}{\Theta(x)}x_1-\frac{25}{\Theta^2(x)}x_2-
\frac{45}{\Theta^3(x)}x_3,\quad \dot x_2=x_1,\quad \dot x_3=x_2.
$$
 Since $\Theta(x)$ is the time motion from point
  $x$ to the origin, that is,
 the equality  $\dot \Theta(x(t))=-1$
  is satisfied, thus,
$\;\Theta(x(t))=\Theta_0-t,\;$ where $\Theta_0$ is a
 positive root of equation (\ref{l2_41}) for $x=x_0.$
 Consequently, the trajectory of the solution is the solution
 of the Cauchy problem of the form
$$ \begin{array}{l}
\displaystyle\dot x_1=-\frac{6}{\Theta_0-t}x_1-
\frac{25}{(\Theta_0-t)^2}x_2-\frac{45}{(\Theta_0-t)^3}x_3,\quad \dot
x_2=x_1,\quad \dot x_3=x_2,\\ x_1(0)=x^0_1,\quad x_2(0)=x^0_2,\quad
x_3(0)=x^0_3.
\end{array} $$
 This equation is reduced to the differential equation of the form
$$\left(\Theta_0-t\right)^3 x_3^{(3)}
+6\left(\Theta_0-t\right)^2\ddot x_3+ 25\left(\Theta_0-t\right)\dot
x_3+45x_3=0,$$
 with initial conditions
 $\; x_3(0)=x^0_3,$
$\; \dot x_3(0)=x^0_2$, and $\ddot x_3(0)=x^0_1.$
 With the change of variables
$\;t=\Theta_0-e^\tau$,
 this Euler-type differential equation is reduced to a differential
 equation with  \ \ constant coefficients \ \  with \ \  respect to
 \newpage
  the function
 $y(\tau)=x_3(\Theta_0-e^\tau),$ which has the form
  $\;y'''-9y''+33y'-45y=0.$
 Thus, we have
$$\;y(\tau)=e^{3\tau}\left(c_1+c_2\cos \sqrt{6}\ \tau
+ c_3\sin\sqrt{6}\ \tau\right),\;$$
 where the constants
  $\;c_1,$ $c_2,$ and $c_3\;$ are found from the conditions
$$y(\tau_0)=x^0_3,\quad y'(\tau_0)=- \Theta_0x^0_2,\quad
y''(\tau_0)-y'(\tau_0)= \Theta^2_0x^0_1,\quad (\tau_0=\ln \Theta_0)
$$
 and are equal  to
$$\begin{array}{l}
\displaystyle c_1=\frac{1}{6 \Theta_0}\left
(\frac{15}{\Theta^2_0}x^0_3+
\frac{5}{\Theta_0}x^0_2+x^0_1 \right),\\[10pt]
\displaystyle c_2=\xi_1 \cos\left(\sqrt{6} \ln \Theta_0\right)-
\xi_2 \sin\left(\sqrt{6} \ln \Theta_0\right),\\[10pt]
\displaystyle c_3=\xi_1 \sin\left(\sqrt{6} \ln \Theta_0\right)+
\xi_2 \cos\left(\sqrt{6} \ln \Theta_0\right).
\end{array}$$
Here
$$\xi_1= -\frac{1}{6 \Theta_0}\left (\frac{9}{\Theta^2_0}\ x^0_3+
\frac{5}{\Theta_0}\ x^0_2+x^0_1\right),\quad
\xi_2=-\frac{1}{\sqrt{6}\  \Theta^2_0} \left(\frac{3}{\Theta_0}\
x^0_3+x^0_2\right).$$
 Since
 $x_3(t)=y(\ln( \Theta_0-t)),$ and the functions $x_2(t)$ and $x_1(t)$
  are found by taking the derivative  of the function
 $x_3(t),$
$$x(t)=\left(\begin{array}{l}
\left( \Theta_0{-}t\right) \left(6c_1+5\sqrt{6}\ \xi_2 \cos\alpha
(t) - 5\sqrt{6}\ \xi_1 \sin\alpha (t)\right)\\[3pt]
\left( \Theta_0{-}t\right)^2\left({-}3c_1{-}\left(3\xi_1{+}
\sqrt{6}\ \xi_2\right)\cos\alpha (t){+}\left({-}3\xi_2{+}\sqrt{6} \
\xi_1\right)
\sin\alpha (t) \right)\\[3pt]
\left( \Theta_0-t\right)^3\left( c_1+\xi_1\cos\alpha (t)+\xi_2
\sin\alpha (t)\right)
\end{array}\right),$$
 where $\;\alpha (t) =\sqrt{6}\ \ln\left( 1-t/ \Theta_0\right).$
  Obviously,  $x(t)\to 0$ for $t\to \Theta_0.$

 The control $u(x) $ on the trajectory $x(t)$ has the form
$$u(x(t))=-6c_1+5(6\xi_1-\sqrt{6}\ \xi_2)\cos\alpha (t)+
5(\sqrt{6}\ \xi_1+6\xi_2) \sin\alpha (t).$$

 For simplicity we solve the problem of
 attaining the origin from the points belonging to the curve
$\;x_1> 0,\;$ $x_2=-41x_1^2/121,$ and $\;x_3=0$.
 In this case, from (\ref{l2_41}), the time of motion
$ \Theta_0 $ from
 $ x_0=(x^0_1,-41(x_1^0)^2/121,0)^* $
is equal to  $41x_1^0/11.$
 The trajectory which begins in this point has the form
$$x(t)=\left(\begin{array}{l}
\displaystyle\frac{41x_1^0-11t}{451}\left(6+5\cos\alpha(t)+
5\sqrt{6}\ \sin\alpha(t)\right)\\[5pt]
\displaystyle
\frac{\left(41x_1^0-11t\right)^2}{9922}\left(-6+4\cos\alpha(t)-
3\sqrt{6}\ \sin\alpha(t)\right)\\[5pt]
\displaystyle
\frac{\left(41x_1^0-11t\right)^3}{327426}\left(6-6\cos\alpha(t)+\sqrt{6}\
\sin\alpha(t)\right)
\end{array}\right),$$
\newpage
where $\alpha(t)=\sqrt{6}\ \ln(1-11t/(41x_1^0)).$
 The control $u(x)$ on this trajectory is determined by the relation
$$ u(t)=-\frac{1}{41}\left(6+
35\cos\alpha(t)\right).$$ Clearly, this control satisfies the
restriction.


\begin{thebibliography}{99}
\bibitem{pontryagin}  {\it Pontryagin, V.G. Boltyanski, R.V. Gamkrelitse  and E.F. Mischenko},
The mathematical theory of optimal processes. Moscow, 1961; English
transls, Gordon and Breach, 1985.

\bibitem{bellman} {\it Bellman, R.E.}   Dynamic Programming. Princeton University
Press, Princeton 1957.

\bibitem {korobov}   {\it Korobov, V.I.}  A general approach to the solution of the bounded control synthesis
problem in a controllability problem, 4(8)  (1979) 582-606.
%Mat. Sb. (N.S.), 1979,  Volume 109(151),    Number 4(8),    Pages
%582–606 (Mi msb2409) This article is cited in 26 scientific papers
%(total in 26 papers)
\bibitem  {korobov_sklyar} {\it Korobov, V.I. and  Sklyar G.M.} Methods for
constructing of positional controls and an admissible maximum
principle. Differentia Equation 26(11) (1990) 1194-1924.

 \bibitem {krein} M.G. Krein, A.A. Nudel'man
\emph{The Markov moment problem and extremal problems}.   \hskip 1em
plus  0.5em minus 0.4em\relax
  Translations of Mathematical Monographs. Vol. 50. 2000.

\end{thebibliography}
\end{document}